\documentclass[a4paper,12pt]{article}

\usepackage[utf8]{inputenc}  
\usepackage{amsmath, amssymb, amsthm,color}  
\usepackage{geometry}  
\usepackage{graphicx}  
\usepackage{hyperref}  
\usepackage{array}
\usepackage{booktabs}
\usepackage{geometry}
\usepackage{longtable}
\usepackage{tabularx,ragged2e}
\usepackage{authblk} 

\newtheorem{theorem}{Teorema}[section]
\newtheorem{lemma}[theorem]{Lema}
\newtheorem{proposition}[theorem]{Proposición}

\newtheorem{definition}[theorem]{Definition}

\theoremstyle{remark}

\newcommand{\F}{\mathbb{F}}

\newcommand{\Z}{\mathbb{Z}}
\newcommand{\R}{\mathbb{R}}

\newcommand{\fl}{\mathrm{fl}}

\begin{document}

\title{Univariate amenable functions}

\author{Carlos Beltrán\thanks{Departamento MATESCO, Universidad de Cantabria. The author is supported by grant PID2020-113887GB-I00 funded by
		MCIN/AEI/ 10.13039/501100011033}}

\date{\today}

\maketitle

\begin{abstract}
The concepts of amenable and compatible functions have been introduced in a recent work, in order to state precise mathematical theorems that guarantee that a backward stable algorithm is also forward stable, and that the composition of two stable algorithms results in an stable algorithm. In this work, we elaborate in this theory for univariate real analytic functions, providing simple tests for both concepts and producing tables for a number of elementary functions which are or fail to be amenable.
\end{abstract}

\section{Introduction}
Numerical analysis has long understood that even algorithms that are mathematically the same can behave very differently when run using floating-point arithmetic. Because of this, the idea of stability was introduced to describe how reliable an algorithm is when working with limited precision. Popular references \cite{Higham2002,TrefethenBau2022} generally distinguish three forms of stability:
\begin{enumerate}
	\item A {\em forward stable} algorithm produces a result that is close to the exact answer for the given input.
	\item A {\em mixed forward-backward stable} algorithm yields an almost exact solution to a slightly perturbed problem.
	\item A {\em backward stable} algorithm computes the exact answer for a nearby problem.
\end{enumerate}
There are clear practical implications of the stability (or the lack of it) for a given algorithm. A classic illustration lies in the computation of matrix eigenvalues. While standard eigensolvers are designed to be stable, the naive approach of computing the characteristic polynomial and subsequently finding its roots is widely regarded as numerically unstable and is strongly discouraged due to its susceptibility to significant errors in floating-point arithmetic (see \cite{BurgisserCuckerCardozo} for a theoretical explanation of this fact).

A fundamental challenge arises when considering the {\em composition of stable algorithms}. Intuitively, one might expect that composing individually stable algorithms would result in a stable overall process. However, this is not always the case; the composition of backward stable algorithms, for instance, can fail to be backward stable, see \cite{Bornemann,BNV}. This observation raises a pertinent and fundamental question: When does composing stable algorithms yield another stable algorithm?

The recent work \cite{BNV} identifies two mild sufficient conditions, termed {\em amenability} and {\em compatibility}, both grounded in the concept of condition numbers, and proves that if two amenable functions $h$ and $g$ can be composed, and they are compatible, then the composition of a forward stable algorithm for $h$ and a forward stable algorithm for $g$ indeed yields a forward stable algorithm for the composite function. Moreover, amenability by itself was also proved to be a relevant concept, since under that hypotheses any mixed forward-backward stable algorithm was proved to be also forward stable, which is a mathematically correct version of the (in)famous rule of thumb 
\begin{quote}
	\begin{center}
	forward error $\lesssim$ backward error $\times$ condition number.
	\end{center}
\end{quote}

The general theory developed in \cite{BNV} was established for functions defined in $\mathbb{R}^n$, and it successfully identified several natural functions (e.g., addition, multiplication of numbers, as well as certain matrix and tensor operations) as amenable. Furthermore, it was shown that while the function $f(x)=\sin(x)$ is not amenable when defined over the entire real line, it does exhibit amenability when restricted to a finite interval with endpoints in $\pi \mathbb{Z}$.

In this work, we provide {\em easily verifiable hypotheses} to guarantee that a univariate real analytic function is amenable. Similarly, we present straightforward conditions to check whether two univariate real analytic functions are compatible. As a direct consequence of these findings, we elaborate on a comprehensive list of amenable and non-amenable univariate functions, offering a practical guide for numerical analysts and practitioners. 

We follow \cite{BNV}, as well as standard definitions, for the concepts and results in the following sections.
\subsection{The floating point arithmetic model}
Recall that a floating point number system $\F\subseteq\R$ with base $2$ is a finite subset of the real numbers of the form $\pm m\cdot 2^{e-t}$, where $t\geq2$ is the precision and $e$ is bounded above and below. The mantissa $m$ is either $0$ or satisfies $2^{t-1}\leq m\leq 2^t-1$. Equivalently, $x\in\F$ if $x=0$ or $x=\pm 2^e\cdot 0.a_1\ldots a_t$ with $a_1\neq0$. For the theory, it is assumed that $e$ is unbounded (that is, $e\in\Z$).

For any $x \in\R$, $\fl(x)$ denotes the number in $\F$ closest to $x$ in absolute value (with the smallest absolute value in case of a tie).

The unit roundoff is $u=2^{-t}$, and $\F$ is denoted by $\F_u$.

With these definitions, $\F_u$ satisfies the following axioms:
\begin{itemize}
	\item For all $x\in\R$, $\fl(x)=x(1+\delta)$ for some $\delta\in\R$, $|\delta|\leq u$.
	\item For $x,y\in \F_u$ and an operation $\circ\in\{+,-,*,/\}$, the floating point operation $x\tilde{\circ}y=\fl(x{\circ}y)$ satisfies
	\[
	x\tilde{\circ}y=(x\circ y)(1+\delta) \text{ for some }\delta\in\R, |\delta|\leq u,
	\]
	with the exception of division by $0$.
\end{itemize}

\subsection{Floating point algorithms}\label{sec:algorithm}
The classical interpretation of an algorithm is a sequence of instructions programmable in languages like Matlab, Python, or C. Formally, an algorithm is a BSS machine, a model of computation that is systematically developed in \cite{BCSS98}, see \cite[Appendix B]{BNV} for a short introduction.

Fix some $S\subseteq\R^m$ and let $f:S\to\R^n$ be a mathematical function. In practice, computing $f(x)$ for a given $x$ requires an algorithm usually performed in a floating point system $\F_u$. Thus, the practical execution of an algorithm, termed a {\em floating point algorithm}, may differ from the intended function and even produce varying outputs on different $\F_u$ or $\F_{u'}$.

For example, the function $f(A) = A^{-1}$ (for invertible matrices $A$) can be computed by different algorithms, such as $\hat f_1$ (using minors and determinant) or $\hat f_2$ (Gaussian elimination with pivoting). While these algorithms should idealistically output the same matrix, namely, the unique inverse of $A$, their floating point implementations yield potentially different results due to the accumulation of rounding errors. Moreover, for a given $\F_u$, the domain of the floating point algorithms $\hat f_{1}^u$ and $\hat f_{2}^u$ may differ from each other and from the domain of $f$. For instance, a singular matrix might yield a non-zero determinant in floating point, allowing $\hat f_{1}^u$ to produce an output. Conversely, a non-defective matrix might conduce to an underflow and thus have a zero determinant in floating point, preventing inversion by $\hat f_{1}^u$. These events may happen for some inputs and some values of $u$ but not for some other inputs or other values, greater or smaller, of $u$.

\begin{definition}[Floating Point Algorithm]\label{def:FPalgorithm}
	Given a BSS machine $\hat f:S\to\R^n$ where $S\subseteq\R^m$, the Floating Point Algorithm $\hat f^{u}$ is the sequence of instructions of $\hat f$, executed in $\F_u$, for any given $u>0$. 
	The initial step of $\hat f^{u}$ involves converting the input to $\F_u$.
\end{definition}

\subsection{The condition number}
For a mathematical function $f:S\to\R^n$, $S\subseteq\R^m$, the condition number $\kappa(f,x)$ of $f$ at $x\in S$ is a measure of the maximal possible variation of the output, given an infinitesimal, arbitrary change in the input. More exactly,
\[
\kappa(f,x)=\lim_{\epsilon\to0}\sup_{y\in S,distance(x,y)\leq \epsilon}\frac{distance(f(x),f(y))}{distance(x,y)}.
\]
We will denote
\[
\mu(f,x)=1+\kappa(f,x).
\]
The condition number is a geometric invariant, independent of the algorithm used to compute $f$, and dependent on the chosen distances in $S$ and the range of $f$. Choosing different distance functions leads to different theories. The distance function that best suits numeric computations done in floating point arithmetic is the relative error distance, which in $\R$ takes the form:
\[
\mathrm{distance}(x,y)=\begin{cases}0&x=y=0\\
	\left|\log\frac{y}{x}\right| &xy>0\\
	\infty&\text{ otherwise}
	\end{cases}
\]
Note that this differs slightly from the classical definition of relative error given by $|x-y|/|x|$. The main advantage of the formula above is that it is an actual mathematical distance, while the relative error is not (since the symmetry condition fails). However, in the case that $x$ and $y$ have the same sign and $\mathrm{distance}(x,y)<<1$ we have 
\[
\mathrm{distance}(x,y)=\left|\log\frac{y}{x}\right|=\left|\log\frac{x+y-x}{x}\right|=\left|\log\left(1+\frac{y-x}{x}\right)\right|\simeq\frac{|x-y|}{|x|},
\]
so the use of the mathematical distance does not perturb our intuition on relative errors. The topology of $\R$ in this metric is the one coming from relative error, so it has three connected components $(-\infty,0)\dot\cup\{0\}\dot\cup(0,\infty)$. 

For a real analytic mapping $f:(a,b)\to\R$, $f\not\equiv0$, with $(a,b)$ a bounded or unbounded open interval, the condition number admits a simple formula
\[
\kappa=\kappa(f,x)=\begin{cases}0&x=0\\
	\infty&x\neq 0 \text{ and }f(x)=0\\
	\frac{|x|\cdot|f'(x)|}{|f(x)|}&\text{ otherwise.}
\end{cases}
\]
The composition law for condition numbers states that if the composition $g\circ h$ is defined, then

\begin{equation}\label{eq:boundmu}
	\kappa(g\circ h,x)\leq \kappa(g,h(x))\kappa(h,x).
\end{equation}
Indeed, if $g$ and $h$ are real analytic with $h(x),g(h(x))\neq0$, by the chain rule we have
\begin{equation}\label{eq:kcomposition}
\kappa(g\circ h,x)=\frac{|x|\cdot|(g\circ h)'(x)|}{|g(h(x))|}=\frac{|h(x)|\cdot|g'(h(x))|}{|g(h(x))|}\frac{|x|\cdot|h'(x)|}{|h(x)|}=\kappa(g,h(x))\kappa(h,x),
\end{equation}
and \eqref{eq:boundmu} becomes an equality.
\subsection{Backward, Mixed and Forward stability}
The three concepts of backward, mixed and forward stability are frequently stated informally (with notations such as $O(\epsilon_{machine})$ that do not specify, for example, if and when the hidden constant can depend on the input). These definitions have been formalized in \cite{BNV} and we specify them below for univariate functions.
\begin{definition}\label{def:bs}
	Let $f:(a,b)\to\R$ with $(a,b)$ a bounded or unbounded open interval, and let $\hat f$ be a BSS machine. Then, $\hat f$ is a {\em backward stable} algorithm for $f$ if there exists a constant $C>0$ such that for all $x\in (a,b)$ we have:
	\[
	0<u<\frac{1}{C}\Rightarrow \exists y\in(a,b):\hat f^{u}(x)=f(y)\text{ and } \mathrm{distance}(x,y)\leq Cu.
	\]
	In particular, an output must be produced for all such choices of $x,u$.
\end{definition}

\begin{definition}\label{def:ms}
	Let $f:(a,b)\to\R$ with $(a,b)$ a bounded or unbounded open interval, and let $\hat f$ be a BSS machine. Then, $\hat f$ is a {\em mixed stable} algorithm for $f$ if there exists a constant $C>0$ such that for all $x\in (a,b)$ we have:
\[
0<u<\frac{1}{C}\Rightarrow \exists y\in(a,b):\mathrm{distance}(\hat f^{u}(x),f(y))\leq Cu\text{ and } \mathrm{distance}(x,y)\leq Cu.
\]
In particular, an output must be produced for all such choices of $x,u$.
\end{definition}

\begin{definition}\label{def:fs}
	Let $f:(a,b)\to\R$ with $(a,b)$ a bounded or unbounded open interval, and let $\hat f$ be a BSS machine. Then, $\hat f$ is a {\em forward stable} algorithm for $f$ if there exists a constant $C>0$ such that for all $x\in (a,b)$ we have:
\[
0<u<\frac{1}{C\mu(f,x)}\Rightarrow \mathrm{distance}(\hat f^{u}(x),f(x))\leq Cu\mu(f,x).
\]
In particular, an output must be produced for all such choices of $x,u$, with a unique exception: if $\mu(f,x)=\infty$ the condition above is vacuous and hence a forward stable algorithm may produce no value, or any value, for an input $x$ such that $\mu(f,x)=\infty$.
\end{definition}

\subsection{Amenability and compatibility of functions}
A property that makes functions amenable for numerical computations with provable stability results has been defined in \cite{BNV}. We write it here for univariate functions.
\begin{definition}
	Let $f:(a,b)\to\R$ with $(a,b)$ a bounded or unbounded open interval. We say that $f$ is amenable if there is a constant $C>0$ such that:
	\begin{itemize}\label{def:amenability}
		\item For all $x\in(a,b)$ with $\kappa(f,x)<\infty$, the interval $B_x=\{y\in\R:\mathrm{distance}(y,x)<1/(C\mu(f,x))\}$ is contained in $(a,b)$. That is, a relative error of size $1/(C\mu(f,x))$ does not move $x$ out of the domain of $f$.
		\item For all $y\in B_x$, we have $\mu(f,y)\leq C\mu(f,x)$. That is, inside $B_x$ the condition number remains essentially bounded by $\mu(f,x)$.
	\end{itemize}
\end{definition}

Also, \cite{BNV} defines the concept of compatibility of functions, designed to give checkable hypotheses for guaranteeing that the composition of stable algorithms yields an stable algorithm. Again, we write the definition here for univariate functions:

\begin{definition}\label{def:compatibility}
	Let $h:(a,b)\to\R$ and $g:(c,d)\to\R$ be real analytic, where $(c,d)$ contains the image of $h$ and both intervals are either bounded or unbounded. We say that $g$ and $h$ are compatible if there exists a constant $C>0$ such that for all $x\in(a,b)$,
	\[
	\mu(g,h(x))\mu(g,x)\leq C\mu(g\circ h,x).
	\]
	That is to say, if a certain reverse inequality of \eqref{eq:boundmu} holds.
\end{definition}

The main outcomes of \cite{BNV} are recalled now:
\begin{theorem}[Prop. 5.6 and Th. 5.7 of \cite{BNV}]\label{theorem: composition of amenable}
	If $g$ and $h$ are compatible amenable functions, then $g\circ h$ is also amenable. Moreover, if $\hat g$, $\hat h$ are respective forward stable algorithms for $g$ and $h$, then $\hat g\circ\hat h$ is a forward stable algorithm for $g\circ h$.
\end{theorem}

\begin{theorem}[Props. 5.10 and 5.11 \cite{BNV}]\label{theorem: chain of implications}
	Let $f$ be amenable. Then:
	\begin{itemize}
		\item Any backward stable algorithm $\hat f$ for $f$ is also mixed stable. 
		\item Any mixed stable algorithm $\hat f$ for $f$ is also forward stable. 
	\end{itemize}
\end{theorem}
The first item of Theorem \ref{theorem: chain of implications} does not require $f$ to be amenable, since it follows immediately from the definitions, but the second item has a nontrivial proof. The results in this paper are designed to allow for an easy use of theorems \ref{theorem: composition of amenable} \ref{theorem: chain of implications} in the case of univariate real analytic functions.

\section{Specific results for one variable functions}
In this section we provide simple forms of the general definitions amenability and compatibility, valid for real analytic, univariate functions defined in real intervals, or unions of real intervals.

\subsection{Sufficient conditions for amenability}
Checking amenability using \cite[Lemma 5.2]{BNV} can be a time consuming task. In this section we show that for the case of real analytic univariate functions a more simple test can be used. The main result is as follows, see Section \ref{sec:proofofpropositions} for a proof.
\begin{proposition}\label{prop:checkamenabilityabbounded}
Let $a,b>0$ and let $f:(a,b)\to\R$ be real analytic and nonzero. Assume that:
\begin{enumerate}
\item Both for $x\to a$ and for $x\to b$ we have $\kappa(f,x)\to\infty$.
\item Both for $x\to a$ and for $x\to b$ we have that $\limsup |H(f,x)|\leq C$ for some constant $C>0$, where
\begin{align}\label{eq:H}
H(f,x)=&\frac{x^2f(x)f''(x)}{f(x)^2+x^2f'(x)^2}.
\end{align}
This claim for $x\to a$ (resp. $x\to b$) is automatically satisfied if $f$ admits an analytic extension to $(a-\epsilon,\infty)$ (resp. $(a,b+\epsilon$)) for some $\epsilon>0$.
\end{enumerate}
Then, $f$ is amenable.
\end{proposition}

In the cases $a=0$ and/or $b=\infty$ the proposition takes a slightly different form (the proof is almost equal and we do not repeat it):
\begin{proposition}\label{prop:checkamenabilityabinfty}
Let $a>0$ and let $f:(a,\infty)\to\R$ be real analytic and nonzero. Let $\mathcal I=\{x\in(a,\infty):f(x)=0\}$. Assume that:
\begin{enumerate}
\item For $x\to a$ we have $\kappa(f,x)\to\infty$, and for any sequence $x_j\to\infty$ such that $\mathrm{dist}(x_j,\mathcal I)\to0$ we have $\kappa(f,x_j)\to\infty$.
\item Both for $x\to a$ and for $x\to \infty$ we have $\limsup |H(f,x)|\leq C$ for some constant $C>0$, where $H$ is given by \eqref{eq:H}. This claim for $x\to a$ is automatically satisfied if $f$ admits an analytic extension to $(a-\epsilon,\infty)$ for some $\epsilon>0$.
\end{enumerate}
Then, $f$ is amenable.
\end{proposition}

\begin{proposition}\label{prop:checkamenabilityazerobbounded}
Let $b>0$ and let $f:(0,b)\to\R$ be real analytic and nonzero. Let $\mathcal I=\{x\in(0,b):f(x)=0\}$. Assume that:
\begin{enumerate}
\item For $x\to b$ we have $\kappa(f,x)\to\infty$, and for any sequence $x_j\to 0$ such that $\mathrm{dist}(x_j,\mathcal I)\to0$ we have $\kappa(f,x_j)\to\infty$.
\item Both for $x\to 0$ and for $x\to b$ we have $\limsup |H(f,x)|\leq C$ for some constant $C>0$, where $H$ is given by \eqref{eq:H}. This claim for $x\to 0$ (resp. $x\to b$) is automatically satisfied if $f$ admits an analytic extension to $(a-\epsilon,\infty)$ (resp. $(a,b+\epsilon)$) for some $\epsilon>0$.
\end{enumerate}
Then, $f$ is amenable.
\end{proposition}

\begin{proposition}\label{prop:checkamenabilityazerobinfty}
Let $f:(0,\infty)\to\R$ be real analytic and nonzero. Let $\mathcal I=\{x\in(0,\infty):f(x)=0\}$. Assume that:
\begin{enumerate}
\item For any sequence $x_j\to 0$ or $x_j\to\infty$ such that $\mathrm{dist}(x_j,\mathcal I)\to0$ we have $\kappa(f,x_j)\to\infty$.
\item Both for $x\to 0$ and for $x\to \infty$ we have $\limsup |H(f,x)|\leq C$ for some constant $C>0$, where $H$ is given by \eqref{eq:H}. This claim for $x\to 0$ is automatically satisfied if $f$ admits an analytic extension to $(-\epsilon,\infty)$ for some $\epsilon>0$.
\end{enumerate}
Then, $f$ is amenable.
\end{proposition}
Finally, if the domain of $f$ is contained in $(-\infty,0)$ the amenability of $f(x)$ is equivalent to the amenability of $f(-x)$ and hence we can apply the propositions above. Moreover, as noted in \cite[Lemma 5.4]{BNV}, if a function is amenable in a finite amount of open intervals, then it is amenable in the union of these intervals. This covers the case of all univariate real analytic functions whose domain has a finite number of connected components.

\subsection{Sufficient conditions for compatibility}
The compatibility condition of Definition \ref{def:compatibility} can be checked in the following simple terms.
\begin{proposition}\label{prop:compatibility}
Let $h:(a,b)\to\R$ and $g:(c,d)\to\R$ be real analytic, where $(c,d)$ contains the image of $h$ and both intervals are either bounded or unbounded. Assume that the following function is bounded in $(a,b)$:
	\[
	\frac{\kappa(h,x)+\kappa(g,h(x))}{1+\kappa(g\circ h,x)}\quad\text{, with the rule }\frac{\infty}{\infty}=1.
	\]
Then $g$ and $h$ are compatible.
\end{proposition}
\section{Some amenable functions}
The main outcome of this section is Table \ref{tab:resultados} where a number of elementary functions are proved to be amenable in their respective domains. The proof consists on a direct application of propositions \ref{prop:checkamenabilityabbounded}, \ref{prop:checkamenabilityabinfty}, \ref{prop:checkamenabilityazerobbounded} and \ref{prop:checkamenabilityazerobinfty}, and is included only for some cases in the following subsections.
\subsection{Univariate polynomials}
Let $f(x)$ be a polynomial. We use Proposition \ref{prop:checkamenabilityazerobinfty}: item 1 is void since $\mathcal I$ is bounded apart from $0$ (which is an isolated point in the topology of relative error) and from $\infty$ (for there are a finite amount of zeros of $f$). As for item 2, note that $H(f,x)$ is the quotient of two polynomials, the degree of the numerator being no greater than the degree of the denominator, and hence $H(f,x)$ is bounded as $x\to\infty$. As for $x\to0$, $H(f,x)$ is bounded from the second item in Proposition \ref{prop:checkamenabilityazerobinfty}. All in one, $f$ is amenable in $(0,\infty)$. The same can be said about the interval $(-\infty,0)$. We conclude that {\em all univariate polynomials are amenable functions in the whole real line}.
\subsection{Rational functions}
Now let $f:\Omega\to\R$, $f(x)=p(x)/q(x)$ with $p,q$ polynomials of respective degrees $d_p$ and $d_q$, and $\Omega=\R\setminus\{x\in\R:q(x)=0\}$. We can consider that the fraction is reduced so that we do not simultaneously have $p(x)=q(x)=0$. We can write down $\Omega$ as a finite union of (possible unbounded) open intervals. If one of such intervals is $(a,b)$ with $a,b>0$, then $q(a)=q(b)=0$. For $x\to a$ we have
\[
f(x)=\frac{p(a)+(x-a)h_1(x)}{C(x-a)^k+(x-a)^{k+1}h_2(x)},
\]
where $k\geq1$, $h_1(x)$ and $h_2(x)$ are polynomials and $C\neq0$. Then,
\begin{multline*}
f'(x)=\frac{h_1(x)+(x-a)h_1'(x)}{(C(x-a)^k+(x-a)^{k+1}h_2(x))}\\-\frac{(p(a)+(x-a)h_1(x))(Ck(x-a)^{k-1}+(k+1)(x-a)^{k}h_2(x)+(x-a)^{k+1}h_2'(x))}{(C(x-a)^k+(x-a)^{k+1}h_2(x))^2}
\end{multline*}
Then, it is immediate to see that in that case
\[
\lim_{x\to a}\kappa(f,x)=\lim_{x\to a}\left|x\frac{f'(x)}{f(x)}\right|=\infty.
\]
The same argument can be applied to $x\to b$, proving that hypotheses 1 of Proposition \ref{prop:checkamenabilityabbounded} holds. Moreover, an elementary but tedious computation shows that:
\[
\lim_{x\to a}H(f,x)=\lim_{x\to b}H(f,x)=\frac{k+1}{k}\leq 2,
\]
and hence it remains bounded. We have thus proved that $f$ is amenable in $(a,b)$.

In an interval of the form $(a,\infty)$ item 1 of Proposition \ref{prop:checkamenabilityabinfty} also holds since $\mathcal I$ is bounded apart from $\infty$. Finally, item 2 in that proposition also holds since $H(f,x)$ is again the quotient of two polynomials, the degree of the numerator being no greater than the degree of the denominator, and hence $H(f,x)$ is bounded as $x\to\infty$. Finally, for $x\to0$, as before, if $p(0)=0$ then $H(f,x)\to (k^2-k)/(k^2+1)$ for some $k\geq1$, and if $q(0)=0$ then we have $f(x)=h(x)/x^k$ for some $k\geq1$ with $h(x)$ a function which is analytic at $0$. Again this implies $H(f,x)\to (k^2+k)/(k^2+1)$ which is bounded by $2$, proving that the hypotheses of \ref{prop:checkamenabilityazerobbounded} also hold. All in one, we have concluded that $f$ is amenable in its whole domain.

\subsection{Algebraic functions}
We now consider $f(x)=x^\alpha$ with $\alpha\neq0$, defined in $(0,\infty)$. Its condition number satisfies
\[
\kappa(f,x)=\begin{cases}0&x=0\\
\alpha&otherwise
\end{cases}
\]
We observe trivially that the hypotheses of Proposition \ref{prop:checkamenabilityazerobinfty} are satisfied and hence $f(x)$ is amenable in $(0,\infty)$.

Now let $f(x)=\sum_{i=1}^ka_ix^{\alpha_i}$ with $a_1,\ldots,a_k\in\R\setminus\{0\}$ and $\alpha_1<\cdots<\alpha_k\in\R$. The domain is $(0,\infty)$  and the condition number is
\[
\kappa(f,x)=\begin{cases}0&x=0\\
\left|\frac{\sum_{i=1}^ka_i\alpha_ix^{\alpha_i}}{\sum_{i=1}^ka_ix^{\alpha_i}}\right|&otherwise
\end{cases}
\]
It is an elementary fact that $f(x)$ has at most a finite number of real zeros, which implies that item 1 in Proposition \ref{prop:checkamenabilityazerobinfty} is void. For item 2, note that
\[
H(f,x)= \frac{\left(a_0+\sum_{i=1}^ka_ix^{\alpha_i}\right)\left(\sum_{i=1}^ka_i\alpha_i(\alpha_i-1)x^{\alpha_i}\right)}{\left(a_0+\sum_{i=1}^ka_ix^{\alpha_i}\right)^2+\left(\sum_{i=1}^ka_i\alpha_ix^{\alpha_i}\right)^2}.
\] 
The limit of this quantity as $x\to0$ is $\alpha_1(\alpha_1-1)/(1+\alpha_1^2)$ and hence $H$ remains bounded as $x\to0$. Similarly, the limit of $H$ as $x\to\infty$ is $\alpha_k(\alpha_k-1)/(1+\alpha_k^2)$ and hence $H$ is also bounded as $x\to\infty$. All in one, item 2 in Proposition \ref{prop:checkamenabilityazerobinfty} also holds and $f$ is amenable in $(0,\infty)$.

\subsection{Trigonometric functions}\label{sec:trigonometric}
We now prove that $f(x)=\sin x$ is amenable when restricted to any interval of the form $(k_1\pi,k_2\pi)$ with $k_1< k_2$ integers. Indeed, this is immediate from Proposition \ref{prop:checkamenabilityabbounded} since
\[
\kappa(\sin,x)=\begin{cases}0&x=0\\
\infty&x=k\pi,0\neq k\in\Z\\
\frac{|x||\cos x|}{|\sin x|}&\text{ otherwise.}
\end{cases}
\]
obviously satisfies $\kappa(\sin,x)\to\infty$ as $x\to k\pi$ with $0\neq k\in\Z$ (and if $k_1$ or $k_2$ are equal to $0$, item 1 of Proposition \ref{prop:checkamenabilityazerobinfty} applies immediately). Item 2 of both propositions also holds for $\sin$ can be analytically extended to $\R$.

However, $f(x)=\sin x$ does not satisfy item 1 of either of the propositions above if the domain contains some unbounded interval, since in that case the sequence $x_j=j\pi+\pi/2$ for $j\to\infty$ satisfies $\mathrm{dist}(x_j,I)=\left|\log\frac{j\pi+\pi/2}{j\pi+\pi}\right|\to0$ as $j\to\infty$, and still $\kappa(\sin,x_j)=0$ for all $j$. Indeed it was proved in \cite[Section 8]{BNV} that $\sin$ is not amenable when defined in an unbounded interval.

A similar claim holds for $f(x)=\cos x$: it is amenable when restricted to any interval of the form $(k_1\pi+\pi/2,k_2\pi+\pi/2)$ with $k_1< k_2$ integers, but it is not amenable if the domain contains some unbounded interval.

We now let $f(x)=\tan x$ defined in $\cup_{k\in\mathcal K}(k\pi-\pi/2,k\pi+\pi/2)$, where $\mathcal K$ is any finite subset of $\mathbb Z$. In order to prove that $\tan$ is amenable, it suffices to see that it is amenable when restricted to any interval of the form $(k\pi-\pi/2,k\pi+\pi/2)$. Hence, we fix any $k\in\mathbb{Z}$ and we check the two items of Proposition \ref{prop:checkamenabilityabbounded}:
\begin{enumerate}
	\item[1.]  The condition number equals
	\[
	\kappa(\tan,x)=\left|\frac{2x}{\sin(2x)}\right|,
	\]
	that goes to $\infty$ as $x\to k\pi-\pi/2$ or $x\to k\pi+\pi/2$.
	\item[2.] An elementary computation shows that:
	\[
	|H(\tan,x)|=\left|\frac{8\,x^2\,\sin^2(4x)}{4\,x^2+\sin^2(2x)}\right|\leq 2.
	\]
\end{enumerate}
We thus conclude that $\tan$ is amenable in  $\cup_{k\in\mathcal K}(k\pi-\pi/2,k\pi+\pi/2)$, $\mathcal K$ being any finite subset of $\mathbb Z$. The same applies to $f(x)=\cot x$ defined in $\cup_{k\in\mathcal K}(k\pi-\pi,k\pi+\pi)$.

\subsection{Exponentials, logarithms and other functions}
In Table \ref{tab:resultados} we have described some elementary functions that satisfy the hypotheses of propositions \ref{prop:checkamenabilityabbounded}, \ref{prop:checkamenabilityabinfty}, \ref{prop:checkamenabilityazerobbounded} or \ref{prop:checkamenabilityazerobinfty} and are thus amenable. In most cases the fact that the hypotheses of the corresponding proposition are satisfied, is a straightforward computation.

An interesting case is that of $f(x)=\sin(\log x)$, which is amenable in $(0,\infty)$. Note that
\[
|H(f,x)|=\left|\frac{\sqrt{2}\,\cos\left(\frac{\pi }{4}+2\,\log\left(x\right)\right)}{2}-\frac{1}{2}\right|\leq \frac{1+\sqrt2}{2},
\]
so one must just check for item 1 of Proposition \ref{prop:checkamenabilityazerobinfty}. Indeed, let $x_j$ be a sequence with $x_j\to0$ or $x_j\to\infty$ and such that $\mathrm{dist}(x_j,\mathcal{I})\to0$. That is to say, we have
\begin{align*}
0=&\lim_{j\to\infty}\mathrm{dist}(x_j,\mathcal{I})\\
=&\lim_{j\to\infty}\min_{k\in\Z}\mathrm{dist}(x_j,e^{k\pi})\\
=&\lim_{j\to\infty}\min_{k\in\Z}\left|\log\frac{x_j}{e^{k\pi}}\right|\\
=&\lim_{j\to\infty}\min_{k\in\Z}\left|\log x_j-k\pi\right|.
\end{align*}
In other words, there exists a sequence $k_j\subseteq \Z$ such that $\epsilon_j=\log x_j-k_j\pi$ satisfies $\epsilon_j\to0$. But then
\begin{align*}
\kappa(f,x_j)=&\left|\frac{\cos\left(\log\left(x_j\right)\right)}{\sin\left(\log\left(x_j\right)\right)}\right|\\
=&\left|\frac{\cos\left(k_j\pi+\epsilon_j\right)}{\sin\left(k_j\pi+\epsilon_j\right)}\right|\\
=&\left|\frac{\cos \epsilon_j }{\sin \epsilon_j }\right|\to\infty.\\
\end{align*}
From Proposition \ref{prop:checkamenabilityazerobinfty}, it follows that $f(x)=\log(\sin x)$ is amenable in $(0,\infty)$.

\section{Some non--amenable functions}
We have seen in Section \ref{sec:trigonometric} that $f(x)=\sin x$ is not amenable in its whole domain $\R$, and indeed that it does not satisfy item 1 of Proposition \ref{prop:checkamenabilityazerobinfty}. We now show an example of a nonamenable function that satisfies item 1 but not item 2 of Proposition \ref{prop:checkamenabilityabinfty}: let
\[
f(x)=1+\sqrt{x-1},\text{ defined in }(1,\infty)
\]
We have
\[
\kappa(f,x)=\frac{x}{2\,\left(\sqrt{x-1}+1\right)\sqrt{x-1}}\to\infty\text{ as }x\to1,
\]
and hence item 1 of Proposition \ref{prop:checkamenabilityabinfty} holds. However,
\[
H(f,x)=-\frac{x^2\,\left(\sqrt{x-1}+1\right)}{\left({x^2}+4{\left(x-1+\sqrt{x-1}\right)}^2\right)\,{\sqrt{x-1}}}
\]
is not bounded as $x\to1$, so item 2 of Proposition \ref{prop:checkamenabilityabinfty} does not hold, and we cannot claim that $f$ is amenable in $(1,\infty)$. Indeed, we now see that $f$ is not amenable in $(1,\infty)$ by checking that there exists no constant $C>0$ with the property that $B_x\subseteq(1,\infty)$ for all $x\in(1,2)$, where
\begin{align*}
	B_x=&\left\{y\in\R:\mathrm{dist}(x,y)\leq\frac{1}{C\mu(f,x)}\right\}\\
	=&\left\{y\in\R:\left|\log\frac{x}{y}\right|\leq\frac{1}{C\mu(f,x)}\right\}\\
	\supseteq&\left\{y\in(0,x):\log\frac{x}{y}\leq\frac{1}{C\mu(f,x)}\right\}\\
	=&\left\{y\in(0,x):y\geq xe^{-\frac{1}{C\mu(f,x)}}\right\}\\
	\ni&\;xe^{-\frac{1}{C\mu(f,x)}}.
\end{align*}
It hence suffices to see that for all $C>6$ there exists $x\in(1,2)$ with $y_x=xe^{-\frac{1}{C\mu(f,x)}}<1$. Fix $C>6$ and note that (using $e^{-t}<1-t/2$ for $t\in(0,1)$):
\begin{align*}
	y_x=&xe^{-\frac{1}{C\mu(f,x)}}\\
	<&x\left(1-\frac{1}{2C\mu(f,x)}\right)\\
	=&x\left(1-\frac{1}{2C\left(1+\kappa(f,x)\right)}\right)\\
	=&x\left(1-\frac{1}{2C\left(1+\frac{x}{2\,\left(\sqrt{x-1}+1\right)\sqrt{x-1}}\right)}\right)\\
	<&x\left(1-\frac{\left(\sqrt{x-1}+1\right)\sqrt{x-1}}{6C}\right)\\
	<&x\left(1-\frac{\sqrt{x-1}}{6C}\right).
\end{align*}
Taking $x=1+C^{-4}$ we see that
\begin{align*}
	y_x<&\left(1+\frac{1}{C^4}\right)\left(1-\frac{C^{-2 }}{6c}\right)\\
	=&\left(1+\frac{1}{C^4}\right)\left(1-\frac{1}{6C^3}\right)\\
	<&\left(1+\frac{1}{C^4}\right)\left(1-\frac{1}{C^4}\right)\\
	=&1-\frac{1}{C^8}<1.
\end{align*}
Hence, the function is not amenable. 

Based on these examples, it is easy to generate functions that fail to satisfy either item 1 or item 2 of any of the propositions \ref{prop:checkamenabilityabbounded}, \ref{prop:checkamenabilityabinfty}, \ref{prop:checkamenabilityazerobbounded} or \ref{prop:checkamenabilityazerobinfty}, see Table \ref{tab:resultados_nonamenable}.

\begin{table}[h!]
	\centering
	\begin{tabularx}{\textwidth}{|c|c|c|c|X|}
		\hline 
		$\sharp$ & $f(x)$ & $\kappa(f,x)$ &  Domain & Failure of hypotheses \\ 
		\hline 
		\rule{0pt}{.7cm} \Centering(1) \label{eq:fila_1} & $\sin\left(x\right)$ & $|x\cot x|$ &$\mathbb{R}$ & Item 1 in Prop. \ref{prop:checkamenabilityazerobinfty}: let $x_k=k\pi+\pi/2$. Then, $\kappa(\sin,x_k)=0$ for all integer $k$ but $\mathrm{dist}(x_k,\mathcal{I})\to0$ as $k\to\infty$ \\ 
		\hline 
		\rule{0pt}{.7cm} \Centering(2) \label{eq:fila_2} & $\cos\left(x\right)$ & $|x\tan x|$ &$\mathbb{R}$ & Item 1 in Prop. \ref{prop:checkamenabilityazerobinfty}: let $x_k=k\pi$. Then, $\kappa(\cos,x_k)=0$ for all integer $k$ but $\mathrm{dist}(x_k,\mathcal{I})\to0$ as $k\to\infty$ \\ 
		\hline 
		\rule{0pt}{.7cm} \Centering(3) \label{eq:fila_3} & $\sqrt{x-1}+1$ & $\left|\frac{x}{2\,\left(x+\sqrt{x-1}-1\right)}\right|$ &$(0,\infty)$ & Item 2 in Prop. \ref{prop:checkamenabilityabinfty}: $\underset{{x\to 1^+}}{\lim}H(f,x)=-\infty$ \\ 
		\hline 
		\rule{0pt}{.7cm} \Centering(4) \label{eq:fila_4} & ${\sin\left(x\right)}^2+1$ & $\left|\frac{2\,x\,\cos\left(x\right)\,\sin\left(x\right)}{{\sin\left(x\right)}^2+1}\right|$ &$(1,\infty)$ & Item 2 in Prop. \ref{prop:checkamenabilityazerobinfty}: $\underset{{x\to \infty}}{\limsup}H(f,x)=\infty$ \\ 
		\hline  
		\hline 
	\end{tabularx} 
	\caption{Some nonamenable functions.} 
	\label{tab:resultados_nonamenable} 
\end{table} 

\section{Proofs of the main results}

\subsection{Proof of Proposition \ref{prop:checkamenabilityabbounded}}\label{sec:proofofpropositions}
We will use Lemma 5.2 of \cite{BNV}, which we reproduce now instantiated for the case of univariate functions:

\begin{lemma}\label{lem:5.2} Given a univariate function $f : (a,b) \to \cup_k \mathbb{R}$ with $0<a<b$, assume that there is a constant $C$ such that the following properties hold:
\begin{itemize}
	\item[(i)] Let $(x_0, x_1, x_2, \ldots) \subseteq (a,b)$ be a sequence such that
	\[
	\text{dist}(x_j, \{a,b\} \cup \mathcal{I}) \to 0,
	\]
	where $\mathcal{I} = \{x \in (a,b) : \kappa(f, x) = \infty \}$ is the ill-posed locus. Then, $\kappa(f, x_j) \to \infty$.
	
	\item[(ii)] Let $V = (a,b) \setminus \mathcal{I}$ be the set where the condition number is finite and let $\mu_f : V \to \mathbb{R}, \, x \mapsto \mu(f, x)$. The condition number of $\mu_f$ (i.e., the level-2 condition number) satisfies
	\begin{equation}\label{eq:kappamu}
		\kappa(\mu_f, x) \leq \frac C4 \mu(f, x), \quad \forall x \in (a,b).
	\end{equation}
\end{itemize}

Then, $f$ is amenable with amenability constant $C$.
\end{lemma}
Moreover, as noted in \cite[Remark 5.3]{BNV}, \eqref{eq:kappamu} is satisfied if
\begin{equation}\label{eq:condicion5.3}
\left|x\frac{d}{dx}\kappa(f,x)\right|\leq C(1+\kappa(f,x))^2
\end{equation}
for some constant $C$. We now prove Proposition \ref{prop:checkamenabilityabbounded}. Let $x_j$ be a sequence such that $x_j\to y\in(a,b)$ and $\mathrm{dist}(x_j,\mathcal I)\to0$, which implies $y\in\mathcal I$. Then, $f(y)=0$ and in some open neighborhood of $y$ we have for some integer $k\geq1$:
\[
f(x)=\sum_{\ell\geq k}(x-y)^\ell a_{\ell},\quad a_k\neq0.
\]
We thus have either $\kappa(f,x_j)=\infty$ (since $x_j\neq0$ and it may happen that $f(x_j)=0$) or otherwise
\begin{align*}
\kappa(f,x_j)=&\frac{|x_j|\cdot|f'(x_j)|}{|f(x_j)|}\\
=&\frac{|x_j|\cdot|\sum_{\ell\geq k}\ell(x_j-y)^{\ell-1} a_{\ell}|}{|\sum_{\ell\geq k}(x_j-y)^\ell a_{\ell}|}\\
=&\frac{|x_j|\cdot|\sum_{\ell\geq k}\ell(x_j-y)^{\ell-k} a_{\ell}|}{|(x_j-y)\sum_{\ell\geq k}(x_j-y)^{\ell-k} a_{\ell}|}\\
\to&\infty,
\end{align*}
the last because the series in the numerator and denominator are convergent in some interval containing $y$. We have proved that for every convergent sequence $x_j\to y\in(a,b)$ and $\mathrm{dist}(x_j,\mathcal I)\to0$, we have $\kappa(f,x_j)\to\infty$. Now, let us prove that any sequence $x_j\in(a,b)$ with $\mathrm{dist}(x_j,\mathcal I)\to0$ satisfies $\kappa(f,x_j)\to\infty$. Indeed, in other case we would have an infinite subsequence $x_{j_k}$ with $\mathrm{dist}(x_{j_k},\mathcal I)\to0$ and $\kappa(f,x_j)\leq C$, some constant $C$. Such a subsequence cannot converge to an interior point $y\in(a,b)$ as we just proved, nor to $a$ or $b$ from the first hypotheses of the proposition. All in one, we have proved that (i) in Lemma \ref{lem:5.2} holds.

As for item (ii) of that result, equivalently to \eqref{eq:condicion5.3}, we must check that the following function is bounded above and below in $(a,b)\setminus\mathcal I$:
\begin{align*}
G(f,x)=&\frac{x\frac{d}{dx}\left(\frac{xf'(x)}{f(x)}\right)}{1+\left(\frac{xf'(x)}{f(x)}\right)^2}\\
=&\frac{xf(x)f'(x)+x^2f(x)f''(x)-x^2f'(x)^2}{f(x)^2+x^2f'(x)^2}.
\end{align*}
Note that  $2|xf(x)f'(x)|\leq f(x)^2+x^2f'(x)^2$. Hence, the first and third terms in the numerator yield bounded expressions. Hence we must only check that $H(f,x)$ is bounded above and below in $(a,b)\setminus\mathcal I$. Since $H(f,x)$ is given by the quotient of two analytic functions, it suffices to check boundedness in the extremes of the interval (which is granted by the second hypotheses of the proposition) and in the inner points where the denominator approaches zero. For this last case, let $x\to y \neq0$ such that $f(y)=f'(y)=0$, so we have for some $k\geq2$:
\[
f(x)=\sum_{\ell\geq k}(x-y)^\ell a_{\ell},\quad a_k\neq0.
\]
Hence,
\begin{align*}
H(f,x)=&\frac{x^2(\sum_{\ell\geq k}(x-y)^\ell a_{\ell})(\sum_{\ell\geq k}\ell(\ell-1)(x-y)^{\ell-2} a_{\ell})}{(\sum_{\ell\geq k}(x-y)^\ell a_{\ell})^2+x^2(\sum_{\ell\geq k}\ell(x-y)^{\ell-1} a_{\ell})^2}\\
=&\frac{x^2k(k-1)a_k^2+h_1(x)}{k^2a_{k}^2x^2+h_2(x)}\\
\to&\frac{k(k-1)}{k^2}\leq1,
\end{align*}
and hence $H(f,x)$ remains bounded in this case. In the argument above, $h_1$ and $h_2$ are two real analytic functions defined in a neighborhood of $y$ and satisfying $h_1(y)=h_2(y)=0$. Note that if $f$ has analytic extension to $(a-\epsilon,a)$ the same argument with $y=a$ shows that the limit of $H(f,x)$ when $x\to a$ exists and is bounded, and similarly for $b$.

This finishes the proof of the result.

\subsection{Proof of Proposition \ref{prop:compatibility}}
From Definition \ref{def:compatibility}, in order to see that $g$ and $h$ are compatible we need to check that
\[
\mu(g,h(x))\mu(h,x)\leq C\mu(g\circ h,x),
\]
for some fixed constant $C\in\R$. That is to say, we need to check that the following is bounded (with the rule $\infty/\infty=1$):
\[
B(x)=\frac{(1+\kappa(g,h(x)))(1+\kappa(h,x))}{1+\kappa(g\circ h,x)}.
\]
Denote by $B(x)$ this last expression. We claim that the boundedness of this function is implied by the boundedness of the more simple expression 
\[
A(x)=\frac{\kappa(g,h(x))+\kappa(h,x)}{1+\kappa(g\circ h,x)}
\]
used in the proposition. Indeed, let us assume that $A(x)\leq C$ is bounded. We distinguish several cases:
\begin{enumerate}
	\item If $x=0$, we have $\kappa(h,0)=\kappa(g\circ h,0)=0$. Hence, $C\geq A(0)=\kappa(g,h(0))$, which implies that $B(0)=1+\kappa(g,h(0))=1+A(0)\leq 1+C$ is also bounded.
	\item If $x\neq0$ and $g(h(x))=0$ we have $\kappa(g\circ h,x)=\infty$, hence $B(x)\leq \infty/\infty=1$ is bounded.
	\item If $x\neq 0$, $h(x)\neq 0$ and $g(h(x))\neq0$, then \eqref{def:compatibility} implies $B(x)=1+A(x)\leq 1+C$.
	\item It only remains to analyze the case $x\neq0$, $h(x)= 0$ and $g(h(x))\neq0$, but this is discarded by hypotheses since we would have $\kappa(g,h(x))=0$, $\kappa(h,x)=\infty$ and $\kappa(g\circ h,x)<\infty$, which contradicts the boundedness of $A(x)$.
	
\end{enumerate}
This finishes the proof of the proposition.

\begin{table}[h!]
\centering
\begin{tabularx}{\textwidth}{|c|c|c|c|c|X|}
\hline 
$\sharp$ & $f(x)$ & $\kappa(f,x)$ & Bound on $H(f,x)$ & Domain & Comments \\ 
\hline 
\rule{0pt}{.7cm} \Centering(1) \label{eq:fila_1} & $p\left(x\right)$ & $\left|\frac{x\,\frac{\partial }{\partial x} p\left(x\right)}{p\left(x\right)}\right|$ & Not shown &$\mathbb{R}$ & $p(x)$ any polynomial \\ 
\hline 
\rule{0pt}{.7cm} \Centering(2) \label{eq:fila_2} & $g\left(x\right)$ & $\left|\frac{x\,\frac{\partial }{\partial x} g\left(x\right)}{g\left(x\right)}\right|$ & Not shown &$\mathbb{R}\setminus\{x:q(x)=0\}$ & $g(x)$ any rational function \\ 
\hline 
\rule{0pt}{.7cm} \Centering(3) \label{eq:fila_3} & $\sum_{i=1}^ka_i x^{\alpha_i}$ & $\left|\frac{x\,\frac{\partial }{\partial x} g\left(x\right)}{g\left(x\right)}\right|$ & Not shown &$(0,\infty)$ & $k\in\Z$, $\alpha_i\geq0$, $a_i\in\R$ \\ 
\hline 
\rule{0pt}{.7cm} \Centering(4) \label{eq:fila_4} & $\sin\left(x\right)$ & $\left|\frac{x\,\cos\left(x\right)}{\sin\left(x\right)}\right|$ & $\begin{matrix}{u^2\pi^2,}\\{u=\max(|k_1|,|k_2|)}\end{matrix}$ &$(k_1\pi,k_2\pi)$ & $k_1\leq k_2\in\Z$ \\ 
\hline 
\rule{0pt}{.7cm} \Centering(5) \label{eq:fila_5} & $\cos\left(x\right)$ & $\left|-\frac{x\,\sin\left(x\right)}{\cos\left(x\right)}\right|$ & $\begin{matrix}{(1+u^2)\pi^2,}\\{u=\max(|k_1|,|k_2|)}\end{matrix}$ &$\left(k_1\pi-\frac{\pi}{2},k_2\pi+\frac{\pi}{2}\right)$ & $k_1\leq k_2\in\Z$ \\ 
\hline 
\rule{0pt}{.7cm} \Centering(6) \label{eq:fila_6} & $\mathrm{tan}\left(x\right)$ & $\left|\frac{2x}{\sin(2x)}\right|$ & $2$ &$\left(k_1\pi-\frac{\pi}{2},k_2\pi+\frac{\pi}{2}\right)$ & $k\in\Z$ \\ 
\hline 
\rule{0pt}{.7cm} \Centering(7) \label{eq:fila_7} & $\frac{1}{\sin\left(x\right)}$ & $\left|-\frac{x\,\cos\left(x\right)}{\sin\left(x\right)}\right|$ & $3(1+|k|)^2\pi^2$ &$(k\pi,k\pi+\pi)$ & $k_1\leq k_2\in\Z$ \\ 
\hline 
\rule{0pt}{.7cm} \Centering(8) \label{eq:fila_8} & $\frac{1}{\cos\left(x\right)}$ & $\left|\frac{x\,\sin\left(x\right)}{\cos\left(x\right)}\right|$ & $3\left(\frac12+|k|\right)^2\pi^2$ &$\left(k\pi-\frac{\pi}{2},k\pi+\frac{\pi}{2}\right)$ & $k\in\Z$ \\ 
\hline 
\rule{0pt}{.7cm} \Centering(9) \label{eq:fila_9} & $\mathrm{cot}\left(x\right)$ & $\left|\frac{2x}{\sin(2x)}\right|$ & $2$ &$(k\pi,k\pi+\pi)$ & $k\in\Z$ \\ 
\hline 
\rule{0pt}{.7cm} \Centering(10) \label{eq:fila_10} & ${\mathrm{e}}^x$ & $\left|x\right|$ & $1$ &$\R$ &  \\ 
\hline 
\rule{0pt}{.7cm} \Centering(11) \label{eq:fila_11} & $\log\left(x\right)$ & $\left|\frac{1}{\log\left(x\right)}\right|$ & $1$ &$(0,\infty)$ &  \\ 
\hline 
\rule{0pt}{.7cm} \Centering(12) \label{eq:fila_12} & $\mathrm{sinh}\left(x\right)$ & $\left|\frac{x\,\mathrm{cosh}\left(x\right)}{\mathrm{sinh}\left(x\right)}\right|$ & $1$ &$(-\infty,\infty)$ &  \\ 
\hline 
\rule{0pt}{.7cm} \Centering(13) \label{eq:fila_13} & $\mathrm{cosh}\left(x\right)$ & $\left|\frac{x\,\mathrm{sinh}\left(x\right)}{\mathrm{cosh}\left(x\right)}\right|$ & $1$ &$(-\infty,\infty)$ &  \\ 
\hline 
\rule{0pt}{.7cm} \Centering(14) \label{eq:fila_14} & $\mathrm{atan}\left(x\right)$ & $\left|\frac{x/(x^2+1)}{\mathrm{atanh}(x)}\right|$ & $1$ &$(-\infty,\infty)$ &  \\ 
\hline 
\rule{0pt}{.7cm} \Centering(15) \label{eq:fila_15} & $\mathrm{asinh}\left(x\right)$ & $\left|\frac{x/\sqrt{x^2+1}}{\mathrm{asinh}(x)}\right|$ & $1$ &$(-\infty,\infty)$ &  \\ 
\hline 
\rule{0pt}{.7cm} \Centering(16) \label{eq:fila_16} & $\mathrm{acosh}\left(x\right)$ & $\left|\frac{x/\sqrt{x^2-1}}{\mathrm{acosh}(x)}\right|$ & $2$ &$(1,\infty)$ &  \\ 
\hline 
\rule{0pt}{.7cm} \Centering(17) \label{eq:fila_17} & $\Gamma \left(x\right)$ & $\left|x\,\psi \left(x\right)\right|$ & $3$ &$(0,\infty)$ &  \\ 
\hline 
\rule{0pt}{.7cm} \Centering(18) \label{eq:fila_18} & $\psi \left(x\right)$ & $\left|\frac{x\,\psi'\left(x\right)}{\psi \left(x\right)}\right|$ & 1 &$(0,\infty)$ &  \\ 
\hline 
\rule{0pt}{.7cm} \Centering(19) \label{eq:fila_19} & $\sin\left(\log\left(x\right)\right)$ & $\left|\frac{\cos\left(\log\left(x\right)\right)}{\sin\left(\log\left(x\right)\right)}\right|$ & Not shown &$(0,\infty)$ &  \\ 
\hline 
\hline 
\end{tabularx} 
\caption{Some amenable functions.} 
\label{tab:resultados} 
\end{table} 


\end{document}